\documentclass[a4paper,12pt]{article}
\usepackage{amsmath, amssymb}
\usepackage{doc, exscale, fontenc, latexsym, syntonly}
\usepackage{xypic}
\usepackage{amsfonts}
\usepackage{amscd}

\author{Marco Mackaay}
\title{\large{\bf A note on the holonomy of connections in twisted bundles}} 
 
\newtheorem{Def}{Definition}[section]
\newtheorem{Lem}[Def]{Lemma}
\newtheorem{Thm}[Def]{Theorem}

\newtheorem{Rem}[Def]{Remark}
\newenvironment{pf}{{\bf Proof:}}{\hfill$\Box$\mbox{}}

\def\H{{\cal H}}
\def\G{{\cal G}}

\def\U{{\cal U}}

\def\P{{\cal P}}
\def\L{{\cal L}}

\def\Co{\mathbb{C}}

\def\PO{\P\exp\int}

\def\gpd{\,\lower3pt\hbox{$\longrightarrow$}\hskip-.267in\raise3pt
             \hbox{$\longrightarrow$}\,}

\date{}

\begin{document}

\vspace*{3.5cm}
\maketitle
\thispagestyle{empty}

\begin{flushright}
\parbox{11.3cm}{
\centerline{\small\bf R\'esum\'e}
Twisted vector bundles with connections have appeared in several 
places (see~\cite{BM00,Kap00} and references therein). In 
this note we consider twisted principal bundles with connections and study 
their holonomy, which turns out to be most naturally formulated in terms of 
functors between categorical groups.

{\noindent\bf  Classification A. M. S. :}  18B40, 18F99, 22A22, 53B99, 53C99.}
\end{flushright}

\section*{Introduction}

Let $M$ be a connected finite-dimensional smooth manifold. 
The holonomy map corresponding to a connection in a principal $G$-bundle 
over $M$ yields a smooth group homomorphism 
$$\H\colon\pi_1^1(M)\to G,$$
where $\pi_1^1(M)$ is the {\it thin fundamental group} of $M$ and 
$G$ a Lie group. In Section~\ref{ce} the reader can find 
the precise definition of thin homotopy, but intuitively a homotopy between 
two loops is thin if it does not sweep out any area. 
This formulation of holonomy is due to 
Barrett~\cite{Ba91} (see also Caetano and Picken's work~\cite{CP94}), who also 
proved that one can reconstruct both the bundle and the connection from the 
holonomy.  

A natural question is whether there are generalizations of Barrett's 
results which involve ``higher thin homotopy types'' of $M$. Caetano and 
Picken~\cite{CP98} defined higher thin homotopy groups of $M$, denoted 
$\pi_n^n(M)$. One way to describe the full homotopy 2-type of $M$, which 
contains more information than just $\pi_1(M)$ and $\pi_2(M)$, 
is by means of the {\it fundamental categorical group} of $M$, denoted 
$C_2(M)$ (see Section~\ref{ce}). 
Recall that a categorical group is a group object in the category of 
groupoids. Picken and the present author~\cite{MP01} 
defined the {\it thin fundamental categorical group}, denoted $C_2^2(M)$, 
which encodes the information about the thin homotopy 2-type of 
$M$. In that same paper we showed that if $M$ is simply connected, 
then a smooth group homomorphism $\pi_2^2(M)\to U(1)$ corresponds precisely 
to the holonomy map of a $U(1)$-gerbe with gerbe-connection. If $M$ is 
not simply-connected we showed that the gerbe-holonomy can be described as 
a smooth functor between $C_2^2(M)$ and a certain categorical 
group over $\pi_1^1(M)$, derived from the canonical line-bundle over 
the loop space of $M$ corresponding to the gerbe. 

This paper is about the next question: given an arbitrary 
categorical Lie group, $\G$, what geometrical structure on $M$ 
yields a holonomy functor between $C_2^2(M)$ and $\G$? 
Theorem~\ref{MThm} shows that, for transitive 
categorical Lie groups, i.e. for those which contain an arrow between 
any two objects, the answer is {\it twisted principal bundles with 
connection}, which we define in Section~\ref{tbac}.   

\section{Categorical groups}
\label{ce}

\begin{Def} 
A {\rm categorical group} is a group object in the 
category of groupoids. This 
means that it is a groupoid with a monoidal structure (a multiplication) 
which satisfies the group laws strictly. A {\rm categorical Lie group} is 
a group object in the category of Lie groupoids, which means that 
the underlying groupoid is a Lie groupoid and that the 
tensor product defines a smooth operation with smooth inverses. 
\end{Def}

\noindent For some general theory about categorical groups we refer to~\cite{BroSpen76}. 
Throughout the paper, let 
$M$ be a connected finite-dimensional smooth manifold and 
let $*$ be a base-point 
in $M$.
Our first example in this paper is the {\it fundamental categorical group} of 
$M$, denoted $C_2(M)$. We want to work with smooth loops and 
homotopies in $M$, but the 
problem is that their 
composites need not be smooth in general. However, there is a subset of 
smooth loops and homotopies whose composites are smooth: 

\begin{Def}\cite{CP94,CP98,MP01} 
A based loop $\ell\colon [0,1]\to M$ is said 
to have a {\it sitting point} at $t_0\in[0,1]$, if there exists an 
$\epsilon>0$ such that $\ell$ is constant on $[t_0-\epsilon,t_0+\epsilon]$. 
We denote the set of all smooth based loops in $M$ with $0$ and $1$ as 
sitting points by $\Omega^{\infty}(M)$.  

Similarly, a based 
homotopy $H\colon [0,1]\times [0,1]\to M$, which I call a {\rm cylinder}, 
has a sitting point $(s_0,t_0)$, 
if there exists an $\epsilon>0$ such that $H$ is constant on the disc with 
centre $(s_0,t_0)$ and radius $\epsilon$ in $[0,1]\times [0,1]$. The set 
of all smooth based homotopies with all points in 
the boundary of $[0,1]\times [0,1]$ being sitting points is denoted by 
$\Omega^{\infty}_2(M)$. 
\end{Def} 

\noindent In order to define $C_2(M)$, we need to introduce 
the notion of {\it thin homotopy}:
\begin{Def}\cite{Ba91,CP94}
Two loops, $\ell$ and $\ell'$, are called thin 
homotopic if there exists a 
homotopy between them whose rank is at most equal to 1 everywhere, which 
is denoted by $\ell\stackrel{1}{\sim}\ell'$.
\end{Def}
\begin{Def}\cite{Ba91,CP94} The {\rm thin fundamental group} of 
$M$, denoted $\pi_1^1(M)$, 
consists of all thin homotopy classes of elements in $\Omega^{\infty}(M)$. 
The group operation is induced by the usual composition of loops. 
\end{Def}

\noindent We can define $C_2(M)$ as follows:
\begin{Def} 
\label{cylinders} 
The objects of $C_2(M)$ are the elements of $\pi_1^1(M)$, which we temporarily 
denote by $[\gamma]$. 

For any $\alpha,\beta,\gamma,\mu\in\Omega^{\infty}(M)$ and 
for any homotopies $G\colon \alpha\to\beta$ and $H\colon\gamma\to\mu$, 
we say that 
$G$ and $H$ are equivalent if there exist thin homotopies $A\colon\alpha\to
\gamma$ and $B\colon\beta\to\mu$ such that 
\begin{equation}
\label{eqrelcyl}
AHB^{-1}\sim G.
\end{equation} 

The morphisms between 
$[\gamma]$ and $[\mu]$ are 
the equivalence classes of 
$$\bigcup_{\alpha,\beta}\left\{H\colon\alpha\to\beta\ \vert\ 
[\alpha]=[\gamma],\ [\beta]=[\mu]\right\}$$ modulo this equivalence relation. 

The usual compositions of loops and homotopies define the structure of a
categorical group on $C_2(M)$, as proved in~\cite{MP01}. 
\end{Def}

\begin{Rem} The usual definition of the fundamental categorical group of 
$M$ yields a {\it weak} monoidal groupoid, because the objects are taken to 
be the loops themselves rather than their thin homotopy classes. 
In~\cite{MP01} Picken and the author defined this strict model. 
\end{Rem}

\noindent 
Similarly we can define the {\it thin fundamental categorical group} of 
$M$, denoted $C_2^2(M)$.
 
\begin{Def}~\cite{CP98,MP01} Two cylinders, $c$ and $c'$, are called {\rm thin homotopic} if 
there exists a homotopy 
between them whose rank is at most equal to 2 everywhere, which is 
denoted by $c\stackrel{2}{\sim}c'$.
\end{Def}

\begin{Def}\cite{MP01} The categorical group $C_2^2(M)$ is 
defined exactly as $C_2(M)$ 
except that the equivalence relation (\ref{eqrelcyl}) is now 
$$AHB^{-1}\stackrel{2}{\sim} G.$$ 
\end{Def}

Next we show how to construct a categorical group from any 
central extension of groups,  
\begin{equation}\label{CE}
1\rightarrow H\rightarrow E\stackrel{\pi}{\rightarrow} G\rightarrow 1.
\end{equation} 
We first construct the underlying groupoid, denoted $E\times E/H\gpd G$. 
This is a well-known construction due to Ehresmann (see~\cite{Mack87} 
for references). 
\begin{Def} The objects of 
\begin{equation}\label{Lie}
E\times E/H\gpd G
\end{equation}
are the elements of $G$, the morphisms are equivalences classes in 
$E\times E/H$, where the action of $H$ is defined by 
$(e_1,e_2)h=(e_1h,e_2h).$ Let us denote 
such an equivalence class by $[e_1,e_2]$, and consider it 
to be a morphism from $\pi(e_1)$ to $\pi(e_2)$. Composition is defined by 
$[e_1,e_2h][e_2,e_3]=[e_1,e_3h]$, where $h\in H$. The identity morphism or 
unit of $g\in G$ is taken to be $1_g=[e,e]$, for any $e\in E$ such that 
$\pi(e)=g$. The inverse of $[e_1,e_2]$ is $[e_2,e_1]$.
\end{Def}

\begin{Lem}
\label{intlaw}
The group operations on $G$ and $E$ induce 
a monoidal structure on (\ref{Lie}). The tensor product on objects is 
simply the group 
operation on $G$. On morphisms the tensor product is defined by 
$[e_1,e_2]\otimes[e_3,e_4]=[e_1e_3,e_2e_4]$. 
Because $G$ and $E$ are groups, this makes (\ref{Lie}) into 
a categorical group.
\end{Lem}
\begin{pf} Since the extension is central, the composition and 
tensor product satisfy the {\it interchange law}, i.e. 
$$([e_1,e_2]\otimes [e_3,e_4])([e'_1,e'_2]\otimes [e'_3,e'_4])=
([e_1,e_2][e'_1,e'_2])\otimes([e_3,e_4][e'_3,e'_4]),$$
whenever both sides of the equation make sense. The other requirements for 
a monoidal structure follow immediately from the group axioms in $G$ and 
$E$.
\end{pf}

\begin{Lem}
If (\ref{CE}) is a central extension of Lie groups, then (\ref{Lie}) yields 
a categorical Lie group.
\end{Lem}
\begin{pf} It is well-known that $E\stackrel{\pi}{\to}G$ defines a 
principal $H$-bundle. See~\cite{Mack87} for a proof that (\ref{Lie}) is a 
locally trivial Lie groupoid for any principal $H$-bundle. Clearly the 
tensor product and the inverses are smooth as well. 
\end{pf}

\noindent Clearly we can recover $E$ from (\ref{Lie}) by 
considering the subgroup of all morphisms of the from $[1,e]$ with the 
tensor product as group operation. The target map then defines the 
projection onto $G$ with kernel $H$. There is a simple characterization of 
categorical groups coming from central extensions. 

\begin{Def} A categorical group is called {\rm transitive}, if there is a 
morphism between any two objects.
\end{Def}

\begin{Lem} There is a bijective correspondence between transitive 
categorical (Lie) groups and central extensions of (Lie) groups. 
\end{Lem}
\begin{pf} An arbitrary transitive categorical group, $\G$,  
corresponds, in the way explained above, to the central extension 
$$\G_1(1,\bullet)\stackrel{t}{\to}\G_0,$$
where $\G_1(1,\bullet)$ is the set of all 1-morphisms starting at the 
unit object, $\G_0$ is the set of all objects and $t$ is the target map. 
The transitivity ensures that $t$ is surjective. In any categorical group 
$t$ is a group homomorphism and the interchange law, mentioned already 
in the proof of Lemma~\ref{intlaw}, 
ensures that $\G_1(1,1)$ is central in $\G_1(1,\bullet)$. 
\end{pf}

\noindent Note that $C_2(M)$ and $C_2^2(M)$ are transitive if and only if 
$M$ is simply connected. 

\begin{Rem}
If $\G$ is not transitive, then it does not correspond to a 
central extension, but to something more general called a 
{\rm crossed module}. For an explanation we refer to \cite{BroSpen76}. 
\end{Rem}

\section{Twisted principal bundles and connections}
\label{tbac}
A {\it twisted bundle} is a geometric structure whose failure to be a bundle 
is defined by an abelian \v{C}ech 2-cocycle. They appear in the 
literature in several places~\cite{BM00,Kap00}. 
In this section I have tried to give 
a systematic exposition of some basic facts about twisted bundles and 
connections, using Brylinski's construction~\cite{Bry93} 
of the abelian gerbe which expresses the obstruction to lifting a principal 
$G$-bundle to a central extension $E$ of $G$. Nothing in 
this section is new strictly speaking, but I hope that writing out 
everything explicitly is useful for the reader.

Let ${\cal U}=\{U_i : i\in{\cal I}\}$ be a {\it good} covering 
of $M$ of open sets, i.e. all intersections 
$$U_{i_1\ldots i_n}=U_{i_1}\cap\cdots\cap U_{i_n}$$ 
of elements of $\U$ are contractible or empty. 
From now on we fix a central extension of Lie groups, denoted as in 
(\ref{CE}).   

\begin{Def}
\label{tb}
A {\rm twisted principal $E$-bundle}, usually denoted $\P$, 
consists of a principal 
$G$-bundle, $P$, and a set of 
local principal $E$-bundles $Q_i\stackrel{q_i}{\to}U_i$, which 
allow for the natural projections $Q_i\stackrel{p_i}{\to}Q_i/H$,  
together with a set of bundle isomorphisms 
$\theta_i\colon Q_i/H\to P_i=P|_{U_i}$ and 
a set of bundle isomorphisms $\phi_{ij}\colon Q_i|_{U_{ij}}\to 
Q_j|_{U_{ij}}$ such that $\phi_{ji}=\phi^{-1}_{ij}$ holds and 
the following diagram commutes:
\begin{equation}
\label{diagtb1}
\begin{CD} 
Q_i@>{\phi_{ij}}>>Q_j\\
@V{p_i}VV@VV{p_j}V\\
Q_i/H@>{\theta_j^{-1}}\theta_i>>Q_j/H
\end{CD}
\end{equation}
Two twisted principal $E$-bundles, denoted 
$\P=(P,Q_i,\theta_i,\phi_{ij})$ and 
$\P'=(P',Q'_i,\theta'_i,\phi'_{ij})$, are {\rm equivalent} if there are bundle 
isomorphisms $\psi\colon P\to P'$ and $\phi_i\colon Q_i\to Q'_i$ 
such that the following diagram commutes: 
\begin{equation}
\label{diagtb2}
\begin{CD} 
Q_i@>{\phi_i}>>Q'_i\\
@V{p_i}VV@VV{p'_i}V\\
Q_i/H@>{{\theta'}_i^{-1}\psi\theta_i}>>Q'_i/H
\end{CD}
\end{equation}
\end{Def}
The following lemma is an easy consequence of our definitions and we leave 
its proof as an exercise. 
\begin{Lem}
\label{cech}
The commutativity of (\ref{diagtb1}) implies that there exists a 
smooth \v{C}ech 2-cocycle on $M$ with values in $H$, given by 
local functions $h_{ijk}\colon U_{ijk}\to H$, such that  
$$\phi_{ki}\phi_{jk}\phi_{ij}(q)=qh_{ijk}(q_i(q)),$$
holds, for any $q\in Q_i|_{U_{ijk}}$. 

The commutativity of (\ref{diagtb2}) implies that there exists a 
\v{C}ech 1-cochain on $M$ with values in $H$, given by local functions 
$h_{ij}\colon U_{ij}\to H$, such that 
$$\phi'_{ij}(q)=\phi_j\phi_{ij}\phi_i^{-1}(q)h_{ij}(q_i(q))$$
holds, for any $q\in Q_i|_{U_{ij}}$. Furthermore, the equation 
$$h'_{ijk}\equiv h_{ijk}h_{ij}h_{jk}h_{ki}$$
holds on $U_{ijk}$.
\end{Lem} 

\begin{Rem} 
\label{obsgerbe}
Brylinski~\cite{Bry93} shows that, given a principal $G$-bundle 
and the central extension, there is a canonical $H$-gerbe associated to them, 
whose equivalence class is represented by $h_{ijk}$ in the previous lemma. 
The $Q_i$ in Def.~\ref{tb} are local trivializations of that gerbe and each 
$\phi_{ij}$ is an isomorphism between two different trivializations. As he 
shows, one can always choose $(Q_i,\phi_{ij})$ which define a twisted 
$E$-bundle and any two choices lead to equivalent twisted $E$-bundles. 
\end{Rem}

\begin{Rem}
Choosing trivializations of all $Q_i$ yields a definition of 
the twisted $E$-bundle in terms of smooth functions 
$e_{ij}\colon U_{ij}\to E$ such that $e_{ji}=e^{-1}_{ij}$ and  
$$e_{ij}e_{jk}e_{ki}\equiv h_{ijk}$$
holds on $U_{ijk}$. Similarly, one can express the equivalence of twisted 
$E$-bundles by smooth functions $e_i\colon U_i\to E$ satisfying 
$$e'_{ij}\equiv e_i^{-1}e_{ij}e_jh_{ij}$$
on $U_{ij}$. 
\end{Rem}

\begin{Def} A twisted principal $E$-bundle $\P=(P,Q_i,\theta_i,\phi_{ij})$ 
is called {\rm flat} if the $h_{ijk}$ in Lemma~\ref{cech} are constant 
functions. Two flat twisted principal 
$E$-bundles are called {\rm flat equivalent} if there exists an equivalence 
$(\psi,\phi_i)$ between them, such that the $h_{ij}$ in Lemma~\ref{cech} are 
constant functions.  
\end{Def}

\begin{Rem} From Brylinski's study~\cite{Bry93} of the obstruction gerbe 
already mentioned we deduce at once that a flat twisted principal 
$E$-bundle is equivalent, in the sense of our Definition~\ref{tb}, 
to an ordinary principal $E$-bundle, but not necessarily equal to one. 
\end{Rem}

\begin{Rem} If $M$ is simply-connected, then any transitive 
Lie algebroid, with fibre $\L(E)$, can be integrated to a flat 
twisted principal $E$-bundle according 
to Mackenzie's results in\cite{Mack87} 
on the obstruction theory for integrating 
transitive Lie algebroids to Lie groupoids. \footnote{I thank 
Mackenzie for making this remark after a talk I gave in Sheffield on twisted 
bundles.}
His results also show that equivalent 
Lie algebroids yield flat equivalent flat twisted principal bundles, at 
least if the choice of central extension is the same (one can always mod 
out $H$ and $E$ by a discrete central subgroup, which makes no difference for 
the corresponding Lie algebras of course). There is a good notion of a 
connection in a transitive Lie 
algebroid~\cite{Mack87} and it seems likely that such a 
connection can be integrated 
to a flat connection in the corresponding flat twisted principal bundle, 
as defined below. In that case the results in this paper would provide a 
notion of holonomy for connections in transitive Lie algebroids, even if 
they cannot be integrated to true principal bundles.  
\end{Rem}

Next let us explain what a connection in a twisted principal $E$-bundle, 
$\P=(P,Q_i,\theta_i,\phi_{ij})$, 
is. Following Chatterjee's terminology for connections in 
gerbes~\cite{Ch98}, we distinguish between $0$- and $1$-connections.   
\begin{Def}
\label{ctb}
A {\rm $0$-connection} in $\P$ consists of a 
$G$-connection, $\omega$, in the principal $G$-bundle $P$ and $E$-connections, 
$\eta_i$, in the local principal $E$-bundles $Q_i$, such that  
\begin{equation}
\label{twistconn1}
\theta_i^*p_i^*(\eta_i)=\omega_i=\omega|_{P_i}
\end{equation}
holds, where ${}^*$ denotes the push-forward for connections. 

Two twisted principal $E$-bundles, $\P$ and 
$\P'$, with $0$-connections, 
$(\omega,\eta_i)$ and 
$(\omega',\eta'_i)$ respectively, 
are {\rm equivalent} if there exists an equivalence 
$$(\psi,\phi_i)\colon 
\P\to \P'$$ 
such that 
\begin{equation}
\label{twistconn2}
\psi^*(\omega)=\omega'.
\end{equation}
\end{Def}

\begin{Rem} It might seem that too many $0$-connections are 
equivalent according to the definition above, but that is because we have 
not yet defined $1$-connections nor the equivalence between twisted principal 
bundles with both $0$- and $1$-connection.  
\end{Rem}

\noindent 
In the following lemma we derive two easy consequences of (\ref{twistconn1}) 
and (\ref{twistconn2}), the proof of which we omit. Note that the adjoint 
action of $E$ on $\L(H)$ is trivial and, therefore, for any $i\in {\cal I}$, 
the associated bundle $Q_i\times \L(H)/\sim$ is canonically isomorphic to 
the trivial bundle $U_i\times \L(H)$. Thus any form on $Q_i$, with values 
in the associated bundle above, that vanishes on vertical vectorfields, 
can be canonically identified with a form on $U_i$ with values in $\L(H)$.  

\begin{Lem}
\label{Blem} 
Equation (\ref{twistconn1}) implies that there exists a 
1-form on each $U_{ij}$ with values in $\L(H)$, denoted $A_{ij}$, such that 
$$\eta_j-\phi_{ij}^*(\eta_i)\equiv A_{ij}$$
holds. Furthermore, we have $A_{ji}=-A_{ij}$ and 
$$A_{ij}+A_{jk}+A_{ki}\equiv -h_{ijk}^{-1}dh_{ijk}$$
on $U_{ijk}$. Using the 
trivializations of Lemma~\ref{cech} we get 
$$A_j-e_{ij}^{-1}A_ie_{ij}-e_{ij}^{-1}de_{ij}\equiv A_{ij}$$
on $U_{ij}$.

Equation (\ref{twistconn2}) implies that there exists 
a 1-form on each $U_i$ with values in $\L(H)$, denoted $B_i$, such that 
$$\eta'_i-\phi_i^*(\eta_i)=B_i$$ 
holds. Furthermore, we have  
$$A'_{ij}\equiv A_{ij}+B_j-B_i-h_{ij}^{-1}dh_{ij}$$
on $U_{ij}$.
Using local trivializations we get 
$$A'_i-e_i^{-1}A_ie_i-e_i^{-1}de_i=B_i.$$
\end{Lem}

\begin{Rem} 
Given a $G$-connection in $P$, Brylinski~\cite{Bry93} constructs a 
{\rm connective structure} for the canonical gerbe mentioned in 
Remark~\ref{obsgerbe}. 
This consists of 
a local $\underline{A}^1_{U_i,\L(H)}$-torsor on each $U_i$ together with 
some data relating these different local torsors. A 0-connection in a 
twisted bundle is nothing but an object in each torsor. Brylinski 
shows that $0$-connections always exist.      
\end{Rem}

\begin{Def} Let $\P$ be a twisted bundle with a $0$-connection $\eta_i$. 
A {\rm $1$-connection} in $(\P,\eta_i)$ consists of $2$-forms $F_i$ on 
$U_i$ with values in $\L(H)$ satisfying 
\begin{equation}
\label{twistcurv1}
F_j-F_i\equiv dA_{ij},
\end{equation}
on $U_{ij}$.  

A {\rm connection} in $\P$ consists of a $0$-connection, $(\eta_i)$, and a 
$1$-connection in $(\P,\eta_i)$. 

Two twisted $E$-bundles, $\P$ and $\P'$, with connections, $(\eta_i,F_i)$ 
and $(\eta'_i,F'_i)$ respectively, are {\rm equivalent} if 
there exists an equivalence 
$$(\psi,\phi_i)\colon (\P,\eta_i)\to (\P',\eta'_i)$$ 
such that 
\begin{equation}
\label{twistcurv2}
F'_i=F_i+dB_i,
\end{equation}  
where $B_i$ was defined in Lemma~\ref{Blem}.
\end{Def}

\begin{Rem} A $1$-connection is what Brylinksi~\cite{Bry93} 
calls a {\rm curving} and he shows that there always exists one 
for a given abelian gerbe with connective structure.
\end{Rem}

\begin{Def}
Let $(\P,\eta_i,F_i)$ be a twisted bundle with connection, then 
the global $3$-form on $M$ with values in $\L(H)$ defined by 
$$G|_{U_{i}}=dF_i$$
is called the {\rm curvature}. 
The connection is called {\it flat} if 
$$G\equiv 0$$ 
on $M$. 
\end{Def} 

\noindent The following theorem follows directly from Brylinski's~\cite{Bry93} 
analogous result for abelian gerbes. 

\begin{Lem} A twisted bundle is 
flat if and only if it admits a flat connection. 
\end{Lem}

\section{Holonomy}
\label{H}     

Let $\G$ be the categorical Lie group associated to a central extension of 
Lie groups (\ref{CE}). 
In this section I first show how the 
holonomy of a connection in a twisted principal 
$E$-bundle assembles nicely into a 
functor
$$\H\colon C_2^2(M)\to\G,$$ 
and then I show that this functor contains all the information about 
the twisted bundle and its connection. 

Throughout this section a principal bundle with connection is always 
considered in terms of local forms 
$\P=(e_{ij},g_{ij},h_{ijk},A_i,D_i,A_{ij},
F_i)$, where $(g_{ij},D_i)$ defines a principal $G$-bundle with connection, 
$(e_{ij},A_i,F_i)$ define the local bundles $Q_i$ and the $0$- and 
$1$-connection in $\P$ and $h_{ijk}$ and $A_{ij}$ were defined in the 
Lemmas~\ref{cech} and \ref{Blem}. 

Given a smooth cylinder 
$$c\colon [0,1]^2\to M,$$ 
such that $c(s,0)=c(s,1)=*$ for all $s\in [0,1]$, choose an 
open covering of the image of $c$ in $M$. Let $V_i=c^{-1}(U_i)$, where 
$U_i$ is an open set in the covering of $c([0,1]^2)$. Next choose a 
rectangular subdivision of $[0,1]^2$ such that each little rectangle 
$R_i$ is contained in at least one open set, 
which for convenience I take to be $V_i$. Denote the 
edge $R_i\cap R_j$ by $E_{ij}$, and the vertex $R_i\cap R_j\cap R_k\cap R_l$ 
by $V_{ijkl}$. Let $\epsilon(c)\in U(1)$ be the following complex number:
\begin{equation}
\begin{array}{lll}
\epsilon(c)&=&{\displaystyle\prod_{\alpha}}
\exp\int_{R_{\alpha}}c^*F_{\alpha}\cdot
{\displaystyle\prod_{\alpha,\beta}}\exp\int_{E_{\alpha\beta}}c^*
A_{\alpha\beta}\\
&&\times{\displaystyle\prod_{\alpha,\beta,\gamma,\delta}}
h_{\alpha\beta\gamma}(c(V_{\alpha\beta\gamma\delta}))
h_{\alpha\delta\gamma}(c(V_{\alpha\beta\gamma\delta}))^{-1}.
\label{eq:gholfor}
\end{array} 
\end{equation}

\begin{figure}
$$\xymatrix@!0{&&&&&&&&&&&\\
&*{F_{\alpha}}&&*{F_{\delta}}&&&&&*{A_{\alpha\delta}}&&&\\
\ar[rrrr]&&&&&+&\ar@{-}'[r][rrr]&*{A_{\alpha\beta}}&\ar'[r][rr]\ar'[u][uu]&
*{A_{\delta\gamma}}&&{+}\\
&*{F_{\beta}}&&*{F_{\gamma}}&&&&&*{A_{\beta\gamma}}&&&\\
&&\ar[uuuu]&&&&&&\ar@{-}'[u][uuu]&&&\\
&&&&&&&&&&&\\
&&&&&&&&&&&\\
\ar[rrrr]&&&&&&&&&&&\\
&&&*{g_{\alpha\beta\gamma}g^{-1}_{\alpha\delta\gamma}}\ar@{~>}[ul]&&&&&&&&\\
&&\ar[uuuu]&&&&&&&&&
}$$
\caption{concrete formula for gerbe-holonomy }
\label{concrete}
\end{figure}

\noindent The last two products are to be taken over the labels of 
contiguous faces in the rectangular subdivision only and in such a way that 
each face, edge and vertex appears only once. The convention for the 
order of the labels is indicated in Fig.~\ref{concrete}, and the 
orientation of the surfaces in that picture is to be taken counterclockwise. 
If $c$ is closed, 
then $\epsilon$ is exactly equal to the gerbe-holonomy as Picken and I 
showed in~\cite{MP01}, so, in that case, its value does not depend on any of 
the choices that were made for its definition. 
In general $\epsilon$ depends on the choices that we made in 
(\ref{eq:gholfor}), 
of course. As a matter of fact, its value only depends on the choice of 
covering of the boundary of $c([0,1]^2)$, because changes in the covering of 
the ``middle'' of $c([0,1]^2)$ do not affect $\epsilon$, which can be shown by 
repeated use of Stokes' theorem. 
Let us give one more ill-defined definition. Let 
$$\ell\colon [0,1]\to M$$ 
be a loop, based at $*$, in $M$. Since $g_{ij},D_i$ is an honest 
principal $G$-bundle with connection, one can 
define in the usual way their holonomy along $\ell$, denoted 
by $\H_{0}(\ell)\in H$. When one tries to do the same 
for $e_{ij},A_i$, the usual formula for the holonomy is not well-defined. 
However, this should not stop us. Let the image of $\ell$ be covered by 
certain $U_i$ again, and choose a subdivision of $[0,1]$ such that each 
subinterval $I_i$ is contained in the inverse image of at least one 
open set, taken to be $V_i$. Let $V_{i,i+1}$ be the vertex 
$I_i\cap I_{i+1}$. Define $\H_{1}(\ell)\in E$ as 
\begin{equation}
\label{eq:1hol}
\H_{1}(\ell)=\prod_{i}\P\mbox{exp}\int_{I_i}\ell^*A_i\cdot e_{i,i+1}
(\ell(V_{ij})).
\end{equation}  
In (\ref{eq:1hol}) $\P\mbox{exp}\int$ means the path-ordered integral, 
which one has to 
use because $E$ is non-abelian in general. For the same reason the 
order in the product is important. We are now ready for the definition of the 
holonomy functor $\H$, which of course has to be independent of all choices.
\begin{Def}
\label{hol}
Let $\ell$ represent a class in $\pi_1^1(M)$ and define  
$$\H(\ell)=\H_0(\ell)\in H.$$
As remarked already this is well-defined. 

Let $c\colon [0,1]^2\to M$ represent a class in $C_2^2(M)([\ell],[\ell'])$. Define 
\begin{equation}
\label{Hs}
\H(s)=[\H_1(\ell),\epsilon(s)\H_1(\ell')]\in E\times E/H.
\end{equation}
The next lemma shows that 
this is well-defined indeed.
\end{Def}

\begin{Lem}
\label{hollem}
The holonomy functor $\H$, as defined in Def.~\ref{hol}, is a well-defined 
functor between categorical Lie groups
$$\H\colon C_2^2(M)\to\G,$$
which is independent of all the choices that we made for its definition.
\end{Lem}
\noindent\begin{pf} 
Showing that $\H$ preserves the categorical Lie group structures is very easy, once it 
has been established that it is well-defined. 
Therefore I only show the latter. 
To prove well-definedness one has to show two things: that  
(\ref{Hs}) does not depend on the choice of covering of $c$, and that 
(\ref{Hs}) is equal for all representatives of the equivalence class of $c$.
Let us first prove the first of these two statements. As far as 
$H_1(\ell)$ and $H_1(\ell')$ are concerned, it is clear that only the choice 
of covering of the boundary of $c$ affects their value. We already argued that 
the same is true for the value of $c$, due to Stokes' theorem. It now 
suffices to see what happens when we introduce an extra vertical line in 
our rectangular subdivision and a new covering of the new (smaller) 
rectangles at the boundary. In Fig.~\ref{chc} one can see such a change.
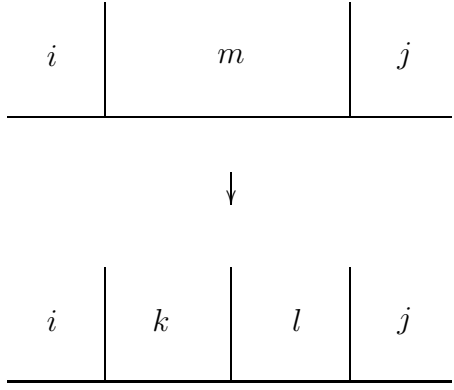
\begin{figure}
$$\xymatrix@-1pc{&&\ar@{-}[dd]+<0pt,-0.25pt>&&&&\ar@{-}[dd]+<0pt,-0.25pt>&&\\
&i&&&m&&&j&\\
\ar@{-}[rrrrrrrr]&&&&&&&&\\
&&&&\ar[d]&&&&\\
&&&&&&&&\\
&&\ar@{-}[dd]+<0pt,-0.25pt>&&\ar@{-}[dd]+<0pt,-0.25pt>&&\ar@{-}[dd]+<0pt,-0.25pt>&&\\
&i&&k&&l&&j&\\
\ar@{-}[rrrrrrrr]&&&&&&&&
}$$
\caption{change in covering}
\label{chc}
\end{figure}

The notation is as indicated in those two pictures. Let us first 
compare the values of $H_1(\ell)$ in the two pictures. In the calculations 
below the pull-back $\ell^*$ has been suppressed to simplify the notation. 
In the first picture 
the part of $\H_1(\ell)$ that matters is equal to  
\begin{equation}
\label{wd1}
e_{im}(V_{im})\cdot\P\exp\int_{I_m}A_m\cdot e_{mj}(V_{mj}),
\end{equation}
and in the second picture that part becomes
\begin{equation}
\label{wd2}
e_{ik}(V_{ik})\cdot\P\exp\int_{I_k}A_k\cdot e_{kl}(V_{kl})\cdot
\P\exp\int_{I_l}A_l\cdot e_{lj}(V_{lj}).
\end{equation}
Now, we can rewrite (\ref{wd2}) to obtain (\ref{wd1}) times an abelian factor. 
In order to do this we have to use the transformation rule for the 
path-ordered integral:
\begin{equation}
\label{po}
\P\exp\int_{[a,b]}(e^{-1}Ae+e^{-1}de)=
e(a)^{-1}\cdot\P\exp\int_{[a,b]}A\cdot e(b).
\end{equation}   
Using the equations satisfied by the $A_i$ as a connection in a 
twisted principal $E$-bundle, we see that (\ref{wd2}) can be rewritten as 
\begin{equation*}
\begin{array}{c}
\displaystyle{
e_{ik}(V_{ik})\cdot\PO_{I_k}(e_{mk}^{-1}A_me_{mk}+e_{mk}^{-1}de_{mk})\cdot 
e_{kl}(V_{kl})}\\[15pt]
\displaystyle{
\times\,\PO_{I_l}(e_{ml}^{-1}A_me_{ml}+e_{ml}^{-1}de_{ml})\cdot 
e_{mj}(V_{mj})\cdot\exp\int_{I_k}A_{mk}\cdot\exp\int_{I_l}A_{ml}.}
\end{array}
\end{equation*}    
Using (\ref{po}) we see that this is equal to
$$
\begin{array}{c}
\displaystyle{
e_{ik}(V_{ik})e_{km}(V_{ik})\cdot\PO_{I_k}A_m\cdot e_{mk}(V_{kl})e_{kl}(V_{kl})
e_{lm}(V_{kl})}\\[15pt]
\displaystyle{\times\,\PO_{I_l}A_m\cdot e_{ml}(V_{mj})e_{lj}(V_{mj})\cdot
\exp\int_{I_k}A_{mk}\cdot\exp\int_{I_l}A_{ml}.}
\end{array}
$$
Finally, using that the coboundary of $e$ is equal to $h$, we get
\begin{equation}
\label{wd3}
\begin{array}{ll}
&\displaystyle{
e_{im}(V_{im})\cdot\PO_{I_k\cup I_l=I_m}A_m\cdot e_{mj}(V_{mj})}\\[15pt]
&\displaystyle{
\times\,h_{ikm}(v_{im})h_{mkl}(V_{kl})h_{mlj}(V_{mj})\cdot\exp\int_{I_k}
A_{mk}\cdot\exp\int_{I_l}A_{ml}}\\[15pt]
=&(8)\times\mbox{abelian factor}.
\end{array}
\end{equation}
Of course a similar calculation can be made for $\H_1(\ell')$. A 
straightforward calculation using Stokes' theorem, which I omit, 
now shows that the inverse of the extra abelian factor in (\ref{wd3}) 
times the extra 
abelian factor in $\H_1(\ell')$ cancel against the extra (abelian) 
factor in $\epsilon(s)$ which appears when it is 
computed for the same change in covering. 

Next let us prove that (\ref{Hs}) is constant on thin homotopy classes. 
Let 
$$c_1,c_2\colon [0,1]^2\to M$$ be two cylinders which represent the 
same class in 
$C_2^2(M)([\ell],[\ell'])$. We have to show that 
\begin{equation}
\label{welldef}
[\H_1(\ell_1),\epsilon(c_1)\H_1(\ell_1')]=
[\H_1(\ell_2),\epsilon(c_2)\H_1(\ell_2')],
\end{equation}
where $\ell_i\stackrel{1}{\sim}\ell$ and $\ell_i'\stackrel{1}{\sim}\ell'$ for 
$i=1,2$. Without loss of generality we may assume that $\ell_2=\ell$ and 
$\ell_2'=\ell'$. Let $A$ be a thin homotopy between $\ell_1$ and $\ell$ and 
let $B$ be a thin homotopy between $\ell_1'$ and $\ell'$, such that 
$$Ac_2B^{-1}\stackrel{2}{\sim}c_1.$$ If the same covering of the boundaries is 
used, then $\epsilon(Ac_2B^{-1})=\epsilon(c_1)$, which follows from the 
general theory of gerbe holonomy developed in~\cite{MP01}. Therefore we have 
$$
\begin{array}{ll}
&[\H_1(\ell_1),\epsilon(c_1)\H_1(\ell_1')]\\[5pt]
=&[\H_1(\ell_1),\epsilon(Ac_2B^{-1})\H_1(\ell_1')]\\[5pt]
=&[\H_1(\ell_1)\epsilon(A)^{-1},\epsilon(c_2)\H_1(\ell_1')\epsilon(B)^{-1}]
\end{array}
$$
It only remains to prove that, for a fixed covering of the boundaries, 
we have 
$$\H_1(\ell_1)\epsilon(A)^{-1}=\H_1(\ell)\quad\mbox{and}\quad 
\H_1(\ell_1')\epsilon(B)^{-1}=\H_1(\ell').$$

\begin{figure}
$$\xymatrix@=0.15cm{*{}\ar@{-}[rrrrrrrr]*{}\ar@{-}[dddddddd]&&*{}
\ar@{-}[dddddddd]
&&*{}\ar@{-}[dddddddd]&&
*{}\ar@{-}[dddddddd]&&*{}\ar@{-}[dddddddd]&&&&*{}
\ar@{-}[rrrrrrrr]*{}\ar@{-}[dddddddd]&&*{}
\ar@{-}[dddddddd]
&&*{}\ar@{-}[dddddddd]&&
*{}\ar@{-}[dddddddd]&&*{}\ar@{-}[dddddddd]&&\\
&&&&&&&&&&&&&&&&&&&&&&\\
*{}\ar@{-}[rrrrrrrr]&&&&&&&&*{}&&&&*{}\ar@{-}[rrrrrrrr]&&&&&&&&*{}&&\\
&&&&&&&&&&&&&&&&&&&&\\
*{}\ar@{-}[rrrrrrrr]&&&&&&&&*{}&\ar[rr]&&&*{}\ar@{-}[rrrrrrrr]&&&&&&&&*{}&
&\mbox{etc.}\\
&&&&&&&&&&&&&&&&&&&&&&\\
*{}\ar@{-}[rrrrrrrr]\ar^2[r]&&*{}\ar^3[d]&&&&*{}&&&&&&*{}\ar@{-}[rrrrrrrr]
\ar^2[r]\ar^8[d]
&*{}&*{}\ar^7[l]\ar^3[r]&*{}&*{}\ar^4[d]&&&&*{}&&\\
*{}&*{}&&&&&&&&&&&*{}&*{}&&&&&&&&&\\
*{}\ar^1[u]\ar@{-}[rrrrrrrr]&&*{}\ar^4[l]&*{}&*{}&*{}&*{}&*{}&*{}&&&&*{}
\ar^1[u]\ar@{-}[rrrrrrrr]&&*{}\ar^6[u]&*{}&*{}\ar^5[l]&*{}&*{}&*{}&*{}&&
}$$
\caption{snake}
\label{sn}
\end{figure}

\noindent Given a rectangular subdivision of $[0,1]^2$ as above, one can 
write the loop around the boundary of $[0,1]^2$ as the composite 
of loops which 
just go around the boundary of one little rectangle $R_i$ at a time 
and are connected 
with the basepoint via a tail lying on some of the edges. 
See Fig.~\ref{sn} for an example, which is like a snake 
(follow the numbers, such that each arrow is numbered on the 
left-hand side). 

Note that the contribution of the 2-forms $F_i$ 
for $\epsilon(A)$ and $\epsilon(B)$ is trivial, because both $A$ and $B$ are 
thin. Using the transformation rule for path-ordered integrals (\ref{po}) 
in the same way as above 
it is not hard to see that $\H_1(\mbox{snake})$ equals  
\begin{equation}
\label{eq1}
\H_1(\ell_1)\epsilon(A)^{-1}\H_1(\ell)^{-1}
\end{equation} 
in the first case and 
\begin{equation}
\label{eq2}
\H_1(\ell_1')\epsilon(B)^{-1}\H_1(\ell')^{-1}
\end{equation}
in the second case. Now recall that on each open set $U_i$ we have an honest 
principal $E$-bundle with an honest connection $A_i$, because 
only globally these data do not match up. Therefore in both cases the 
value of $\H_1$ around the 
boundary of each $R_i$ equals $1$, because $A$ and $B$ are thin. The 
conclusion is that the expressions in (\ref{eq1}) and (\ref{eq2}) are 
both equal to 1 as well.  
\end{pf}

\noindent Clearly the connection in $\P$ is flat if and only if 
$\H$ is constant on ordinary homotopy 
classes of cylinders, which happens if and only if $\H$ defines a 
functor $\H\colon C_2(M)\to\G$. 

Note that the lemma above implies that the element 
$$\H_1(\ell)(\H_1(\ell')\epsilon(c))^{-1}\in E$$ 
is well-defined for any $c\colon [0,1]^2\to M$. 
Kapustin~\cite{Kap00} studied the 
special case in which $E=\mbox{GL}(n,\Co)$ and $c(0,t)=*$ equals 
the trivial loop at the basepoint. His main mathematical result about 
the holonomy of connections in twisted vector bundles seems to 
be that $\mbox{tr}(\H_1(\ell')\epsilon(c))$ is a well-defined complex number 
in that particular case.

In order to understand the sequel, one should note that $G$ acts by 
conjugation both on itself and on $E$. 
\begin{Def}
Given a holonomy functor, $\H$, one can define the conjugate holonomy functor 
$\H^g=g^{-1}\H g$, for any $g\in G$, by 
$$\H^g_0(\ell)=g^{-1}\H_0(\ell)g$$ and 
$$\H^g(c)=[g^{-1}\H_1(\ell_1)g,\epsilon(c)g^{-1}\H_1(\ell_2)g].$$
\end{Def}
As a matter of fact there is a natural isomorphism between $\H$ and $\H^g$ 
defined by 
$$[\H_1(\ell),g^{-1}\H_1(\ell)g]\colon \H_0(\ell)\to \H^g_0(\ell),$$
for any loop $\ell$. Note that this natural isomorphism is well-defined 
indeed. 
\begin{Lem} 
\label{eqh}
Equivalent twisted principal $E$-bundles with connection 
give rise to conjugate holonomy functors. 
\end{Lem}
\noindent\begin{pf} Recall that two twisted principal $E$-bundles with connections, denoted $g_{ij},e_{ij},h_{ijk},A_i,D_i,A_{ij},F_i$ and 
$g'_{ij},e'_{ij},h'_{ijk},A'_i,D'_i,A'_{ij},F'_i$ respectively, 
are equivalent if 
$$
\begin{array}{lll}
g'_{ij}&=&g_i^{-1}g_{ij}g_j\\[5pt]
e'_{ij}&=&e_i^{-1}e_{ij}e_jh_{ij}\\[5pt]
D'_i&=&g_i^{-1}D_ig_i+g_i^{-1}dg_i\\[5pt]
A'_i&=&e_i^{-1}A_ie_i+B_i+e_i^{-1}de_i\\[5pt]
F'_i&=&F_i+dB_i,
\end{array}  
$$
where $p_i\circ e_i=g_i$. 
Thus we see that 
$$\H'_0(\ell)=(\pi g_0(*))^{-1}\H_0(\ell)\,\pi g_0(*),$$
holds, for any loop $\ell$, where by convention $U_0$ is the open set 
that covers the basepoint.

We now have to show that the identity
$$\H'(c)=g_0(*)^{-1}\H(c)g_0(*)$$
holds, for any cylinder $c\colon\ell_1\to\ell_2$. Using the transformation 
rule (\ref{po}) again we get 
$$
\begin{array}{lll}
\displaystyle{\H'_1(\ell_1)}&=&\displaystyle{\prod_{i}\P\mbox{exp}\int_{I_i}
\ell_1^*A'_i\cdot 
e'_{i,i+1}(\ell_1(V_{ij}))}\\[15pt]
&=&\displaystyle{e_0(*)^{-1}\cdot\prod_{i}\P\mbox{exp}\int_{I_i}
\ell_1^*A_i\cdot e_{i,i+1}(\ell_1(V_{ij}))\cdot e_0(*)}\\[15pt]
&&\displaystyle{\times\prod_{i}\mbox{exp}\int_{I_i}\ell_1^*B_i\cdot 
h_{i,i+1}(\ell_1(V_{ij}))}\\[15pt]
&=&\displaystyle{e_0(*)^{-1}\H_1(\ell_1)e_0(*)\times\mbox{abelian factor}}.
\end{array}
$$
Analogously we get 
$$\H'_1(\ell_2)=e_0(*)^{-1}\H_1(\ell_2)e_0(*)\times\mbox{abelian factor}.$$
Using Stokes' theorem it is now easy to see that $\epsilon'(c)$ cancels 
these two abelian factors so that
$$
\begin{array}{lll}
\H'(c)&=&[\H'_1(\ell_1),\epsilon'(c)\H'_1(\ell_2)]\\[5pt]
&=&[e_0(*)^{-1}\H_1(\ell_1)e_0(*),\epsilon(c)e_0(*)^{-1}\H_1(\ell_2)e_0(*)]
\\[5pt]
&=&e_0(*)^{-1}\H(s)e_0(*).
\end{array}
$$  
Finally, the result follows from the observation that conjugation by $g_0(*)$ 
is equal to conjugation by $e_0(*)$, because $E$ is a {\it central} extension. 
\end{pf}

Putting together Barrett's~\cite{Ba91} results about the reconstruction of 
bundles with connections from their holonomies and Picken and my~\cite{MP01} 
analogous results for gerbes, one now almost immediately gets the following 
lemma.

\begin{Lem}
\label{reclem}
Given a smooth functor of categorical Lie groups 
$$\H\colon C_2^2(M)\to \G,$$ 
there exists a twisted principal $E$-bundle with connection 
whose holonomy functor is equal to $\H$. 

If $\H$ is the holonomy functor of a given twisted principal $E$-bundle 
with connection, then the twisted principal 
$E$-bundle with connection constructed from $\H$ is equivalent to the 
given one. 
\end{Lem} 
\begin{pf} Barrett's results~\cite{Ba91} allow us to construct $g_{ij}$ and 
$D_i$ for a given holonomy functor. The rest of the proof relies on the same 
techniques as employed in~\cite{MP01}. We only sketch the construction here. 
Let us see what 
$e_{ij}$ and $A_i$ are in terms of the 
holonomy functor $\H$. In~\cite{MP01} we chose a fixed point 
in each open set, called $x_i\in U_i$, and a fixed point in each double 
overlap, $x_{ij}\in U_{ij}$. We picked a path from the basepoint $*\in M$ to 
each $x_i$ and showed how to fix paths in $U_i$ from $x_i$ to any other 
point in $U_i$ and from $x_{ij}$ to any other point in $U_{ij}$. 
We also showed how to fix homotopies inside $U_i$ 
between any two homotopic paths in $U_i$. In particular we got a fixed 
loop 
\begin{equation*}
\label{loop1}
*\rightarrow x_i\rightarrow y\rightarrow x_j\rightarrow *,
\end{equation*}
for any point $y\in U_{ij}$. Call this loop $\ell_{ij}(y)$. We also got a 
fixed 
homotopy, $c_{ij}(y)$, between $\ell_{ij}(x_{ij})$ and $\ell_{ij}(y)$. 
Now consider 
\begin{equation}
\label{f}
\H(c_{ij}(y))\in E\times E/H.
\end{equation}
Take 
a representative of $\H(c_{ij}(x_{ij}))$ in $E\times E$, which of course is of 
the form $(e,e)$.  
For any $y\in U_{ij}$, take the unique representative of (\ref{f}) of 
the form $(e,e')$ and define  
$$e_{ij}(y)=e'\in E.$$ 
The 
choice which this reconstruction of $e_{ij}$ involves, corresponds to 
gauge fixing. By convention we fix the representative of $\H(c_{ji}(x_{ji}))$ 
to be $(e^{-1},e^{-1})$. Then the identity  
$$e_{ji}(y)=e_{ij}(y)^{-1}$$
follows immediately. Because we have 
$$\ell_{ij}(y)\star\ell_{jk}(y)\star\ell_{ki}(y)\stackrel{1}{\sim}c_*,$$
we can define the 2-cocycle 
$$h_{ijk}(y)=e_{ij}(y)e_{jk}(y)e_{ki}(y)\in\ker\pi\cong H.$$ 
Note that this definition of the 2-cocycle is equal to the one given 
in~\cite{MP01}. Different choices of representatives of 
$\H(s_{ij}(x_{ij}))$ in $E\times E$ yield an equivalent twisted 
principal $E$-bundle.  

Analogously we can reconstruct the 
$A_i$. Given a vector $v\in T_y(U_i)$, we can represent it by a small 
path $q(t)$ in $U_i$, whose derivative at $t=0$ is equal to $v$. Then 
there is a loop
\begin{equation}
\label{ell}
*\rightarrow x_i\rightarrow y\stackrel{q}{\rightarrow} 
q(t)\rightarrow x_i\rightarrow *.
\end{equation}
Call it $\ell_i(q(t))$. For each value of $t$, we can use the fixed homotopy 
in $U_i$ to get a homotopy, $c_i(q(t))$, 
from the trivial loop at $*$ to $\ell_i(q(t))$. 
Consider 
\begin{equation}
\label{A}
\H(c_i(q(t)))=[e(t),e'(t)]\in E\times E/H.
\end{equation} 
Define 
$$A_i(v)=\frac{d}{dt}e(t)^{-1}e'(t)\vert_{t=0}.$$
It is not hard to see that 
$$A_j(v)-e_{ij}(y)^{-1}A_i(v)e_{ij}(y)-e_{ij}^{-1}de_{ij}(v)=A_{ij}(v)$$ 
is exactly the abelian $1$-form that we reconstructed 
in~\cite{MP01}. The fact that the definition of $A_i(v)$ does not depend on 
the particular choice of $q(t)$ follows precisely from the same 
arguments that we used in that paper. 

The reconstruction of the $1$-connection $F_i$ is exactly the 
same as in~\cite{MP01}, because it only depends on the value of the holonomy 
functor around closed cylinders.

The rest of the proof is similar to the proofs of the analogous results for 
bundles and gerbes in~\cite{Ba91} and \cite{MP01} and we omit the details.
\end{pf}

\noindent Choosing a different basepoint in $M$ and a path from that 
basepoint to $*$ yields an equivalence between the two respective 
thin fundamental categorical groups. 
This equivalence induces an equivalence relation on those holonomy 
functors $\H$ which correspond to the same equivalence class of 
twisted principal $E$-bundle with connection. Just as for ordinary 
connections, two holonomy functors are equivalent if and only if they are 
conjugate by an element in $G$.  Together 
with Lem.~\ref{hollem}, Lem.~\ref{eqh} and Lem.~\ref{reclem} these remarks 
prove the following theorem:

\begin{Thm}
\label{MThm}
There is a bijective correspondence:
$$\left\{\mbox{twisted principal $E$-bundles on $M$ with connection}\right\}/
\,\sim $$
$$\longleftrightarrow $$
$$\mbox{Hom}(C_2^2(M),\G)/\G_0.$$
\end{Thm}

\vskip0.5cm
\centerline{\bf Acknowledgements}
\vskip.5cm
\noindent 
I thank Louis Crane and Nigel Hitchin for their, independent, suggestions to 
look at twisted bundles and twisted $K$-theory respectively.
 
I am presently on leave from the Universidade do Algarve 
and working with a postdoctoral fellowship from FCT at the University of 
Nottingham (UK). I thank the former for their financial support and the 
latter for their hospitality. 

This work was also supported by {\em Programa Operacional
``Ci\^{e}ncia, Tecnologia, Inova\c{c}\~{a}o''} (POCTI) of the
{\em Funda\c{c}\~{a}o para a Ci\^{e}ncia e a Tecnologia} (FCT),
cofinanced by the European Community fund FEDER.

\end{document}